\documentstyle{amsppt}

\magnification=1200 \NoBlackBoxes \hsize=11.5cm \vsize=18.0cm

\def\inv{^{-1}}

\def\X{\Cal X}

\def\Y{\Cal Y}

\def\C{\Bbb C}

\def\P{\Bbb P}

\def\ss{\vskip.15in}

\def\L{\Cal L}
\def\O{\Cal O}

\def\Sym{\text{Sym}}

\def\rk{\text{rk}}

\def\p{\partial}

\def\1/2{\frac{1}{2}}

\def\SD{\Cal S^{\text{v}}}

\def\Y{\Cal Y}
\def\SS{\Cal S}

\def\im{\text{im}}

\def\2{{[2]}}
\def\l{\ell}
\topmatter
\title On the geometric genus of subvarieties
of generic hypersurfaces\endtitle
\author Herbert Clemens\footnote
{\raggedright{Partially supported by NSF grant
DMS-9970412}}and
Ziv Ran\footnote{\raggedright{
Partially supported by NSA Grant MDA904-02-1-0094} }
\endauthor

\date Sunday, April 14, 2002 at 20:58\enddate
\address University of Utah\endaddress
\email clemens\@math.utah.edu\endemail
\address University of California, Riverside\endaddress
\email ziv\@math.ucr.edu\endemail
\rightheadtext {Geometric Genus}
\leftheadtext{}
\abstract We prove some
 lower bounds on  certain twists
 of the canonical bundle of a
subvariety of a generic hypersurface in projective
space. In particular we prove that the generic sextic
threefold contains no rational or elliptic curves
and no nondegenerate curves of genus 2.\endabstract

\endtopmatter\document
The geometry of a desingularization
$$Y$$ of an arbitrary
$k-$dimensional subvariety of a generic hypersurface
$$X$$ in an ambient variety
$$W \ \ (\text{e.g.}\ W = \P^{n})$$ has
received much attention over the past decade or so. The first
author ([C], see also [CKM], Lecture 21, for an exposition and
amplification) has proved that for $$k =1 , n =4,
 W = \P^4$$ and $X$ of degree $$d\geq 7,$$
 $Y$ has genus $$g \geq
1,$$ and conjectured that the
same is true for $d=6$ as well;
this conjecture is sharp in the sense that hypersurfaces
of degree $d\leq 5$
do contain rational curves.
What the statements mean in plain terms is that,
in the indicated range, a nontrivial function-field solution
of a generic polynomial equation
must have genus $>0.$ We shall
refer to this as the sextic conjecture,
to distinguish it from
his other conjecture concerning
rational curves on quintics.\par
The sextic conjecture
was proved by Voisin [VE] who showed
more generally
that, for $X$ of degree
$d$ in
$$\P^{n}, n\geq 4, k\leq n-3,$$
$$p_g(Y)>0\ \ \text{if}\ \ d\geq 2n-1-k$$ and $K_Y$
separates generic points if
$$d\geq 2n-k.$$
In the case of codimension 1 in $X$, i.e. $k=n-2$,  Xu
[X] gave essentially sharp geometric genus bounds.
For $X$ a
generic complete intersection of type $(d_1,...,d_k)$
in any
smooth polarized $(n+k)$-fold $M$, Ein [E] proved that
$$p_g(Y)>0$$
if
$$d_1+...+d_k\geq 2n+k-m-1$$
and $Y$ is of
general type if
$$d_1+...+d_k\geq 2n+k-m+2.$$ Ein's bounds are
generally not sharp, e.g. they fail to yield
the sextic conjecture for
$M=\P^{n}$. The paper [CLR] gives some refinements and
generalizations of the results of Ein and Xu by a method which
seems to yield essentially sharp bounds in codimension 1 but not
necessarily in general.
\par
In this paper we give a result which improves
and generalizes
the sextic conjecture.
It gives a lower bound on a twist of the canonical bundle
of an arbitrary subvariety
of a generic hypersurface  in projective
space. For the case of the sextic in $\P^4$
it shows that the minimal genus of a curve is at least 2,
and at least 3 if the curve is nondegenerate. We proceed
to state the result.
\par
First, we fix a projective space $$\P=\P^{n}.$$
We denote by $$\L_d$$ the space of homogeneous polynomials
of degree $d$ on $\P$.
If $Y$ is a variety with a given morphism
$$f:Y\to\P$$ and $A$ is a coherent sheaf on $Y$, the sheaf
$$A\otimes f^*(\O_\P(i))$$ will be denoted by $A(i).$
\proclaim{Theorem 0.1}
Let $$X\in\L_d$$ be generic with
$$d(d+1)/2\geq 3n-1-k, d\geq n, $$
 and
$$f:Y\to X$$ a
desingularization of an irreducible subvariety
of dimension $k$.
Set
$$t=\max(0,-d+n+1+[\frac{n-k}{2}]).$$
Then  either
$$h^0(\omega_Y(t))> 0\tag 0.1$$
or $f(Y)$ is contained in the union of the lines
lying on $X$.\par
In the case where $k=n-3$ we have furthermore:\par
(i)if $h^0(\omega_Y(t))=0$,
then $f(Y)$ is ruled by lines;\par
(ii) if $ d\geq n+2$
and $f(Y)$ is not ruled by lines, then
$$h^0(\omega_Y)+1
\geq\min(\dim({\text{\rm span}}( f(Y))),4).\tag 0.2$$
Finally in the case in which $Y$ is a
curve on a generic sextic threefold,
we  have
$$g(Y)\geq 3,$$
 except possibly if $Y$ is a genus-2 curve
such that $f(Y)$ spans a hyperplane.
\endproclaim\ss
\remark{Remark 0.2}
In case either $k=n-2$
or $k=n-3,\ d>n+2$,
our statement is weaker than the
result of [CLR]. Therefore from now on we may assume
that $k\leq n-3,$ and that if equality holds then
$d=n+2.$\endremark

\remark{Remark 0.3} If $k=2n-2-d\geq 1,\ d\geq n+3,$ our
result implies that $H^0(\omega_Y)\neq 0$ unless $Y$
is ruled by
lines (which in this case means $Y$
is a component of the union
of lines in $X$). In this form the result was first
obtained by
Pacienza [P].
In fact in this paper we build on the
methods of Voisin [V] and Pacienza [P],
which play a very important role in
our results.
If $k>2n-2-d\geq 0$, then of course $Y$ cannot be ruled
by lines, so we conclude in this case that
(0.1) (or, if $k=n-3, d>n+1$, (0.2)) holds.
\endremark
Our proof is based in the observation, already made
by Voisin and, in a more limited context,
by the second author, that the failure of the estimate
(0.1) implies the existence of a line
$$\l=\l(y)$$
through a general point
$$y\in Y=Y_F$$ with the property
that the space of first-order variations
$$F'\ \text{ of}\ F$$
such that
$$Y_{F'}\ni y\ \text{and}\ F'|_\l=F|_\l$$
is larger than expected. This leads us to consider
a distribution $$T'$$ on the space $Y_x$ of
polynomials $F$ such that $Y_F\ni x$. Inspired
by Voisin,  we show in \S 6
that
$T'$ is integrable (Lemmas 6.1-6.2) .
 This in turn leads to the conclusion
that $\l(y)$ is contained in the hypersurface $F$
(Lemma 6.4). We then show in \S 7 that
if $k=n-3$ or, more generally, if a certain
infinitesimal invariant $s'$ vanishes, then in fact, $\l(y)$
is contained in $Y_F$ itself, so that $Y_F$ is ruled by
lines.\par
The main argument is presented in \S 6, 7 which are preceded
by 5 preliminary sections. In \S 1 we give precise
definitions and statements. In \S 2 we recall some
(probably well known) global generation results
for some standard homogeneous vector bundles on
projective spaces and Grassmannians. In \S 3 we
study quotients of the space of homogeneous
polynomials of degree $d$ vanishing at a fixed
point $x$,
and in particular the question of how much of
such a quotient is generated by a generic collection of
degree-$(d-1)$ polynomials. It turns out that if
the answer is less than expected, this can sometimes
by 'explained' by the existence of some special
line $\l$ through $x$. In \S 4 we use the results
of \S 3 to reduce the proof of the main theorem,
in effect, to the study of the distribution $T'$.
\S 5 recalls, and reproves, some well-known results
on canonical bundles of varieties parametrizing
osculating lines to hypersurfaces.\par
This paper supersedes an earlier eprint by the second author
entitled 'Beyond a conjecture of Clemens' (math.AG/9911161).

\heading\bf{ 1. The formal setting}\endheading
We will fix integers
$d,n$ and denote by $\L$ or
$\L_d$ the space of homogeneous forms
of degree $d$ on $\P^n$ and set
$$S=\L-\{0\}.$$ An element $F\in\L$
is written in the form
$$F=\sum\limits_{|A|=d}a_AX^A.$$
We denote by
$$\X\subset S\times\P^n$$
the universal hypersurface,
which is smooth, as is the natural map
$$s:\X\to\P^n,$$
 whose fibre over $x\in\P^n$ we denote
by $X_x,$
while the fibre over $F\in S$ is denote $X_F.$\par
We will be
studying a certain kind of
versal families of $k$ folds on $\X/S$.
By this we
mean the following. Let $T$ be an irreducible
subvariety of
the Chow variety (or Hilbert scheme) of $\P^n$
whose general element corresponds to an irreducible
$k-$dimensional variety. Suppose that
$$\psi:T\to S$$
be a dominant morphism with the property
that if we denote the universal cycle (or subscheme)
over $T$ by $\Cal Z_T,$ then
$$\Cal Z_T\subset\X_T:=\X\times_ST.\tag 1.1$$
Clearly there is no loss of generality
by assuming $T$ is saturated with respect to
projective equivalence, i.e. is stable under
the natural action of GL$(n+1)$.
Note that GL$(n+1)$ acts on (1.1) and choose a
GL$(n+1)$-stable subvariety $T'\subset T$
so that the map $T'\to S$ is \'etale. Then choose
a GL$(n+1)$-equivariant resolution of
$$\Cal Z_{T'}\to T'$$
and denote this by $\Y/S'.$
Replacing $S'$ by a GL$(n+1)$-stable open subset,
we may assume $\Y/S'$ is
a smooth projective relatively $k-$dimensional
family endowed with an $S'-$morphism
$$f:\Y\to \X\times_SS'.$$
We will call this a {\it{versal family of $k-$folds.}}

\heading 2.
\bf{Positivity results}\endheading
\noindent (2.0)\ \  A useful remark is the following.
Let
$$0\to A\to B\to C\to 0$$
be a short exact sequence of vector bundles
of respective ranks $a,b,c.$ Then for any $k>0,$
the exterior power $$\bigwedge\limits^k B$$
is endowed with
a descending filtration $F^.$   with
$$F^i/F^{i+1}\simeq \bigwedge^{k-i}C\otimes\bigwedge
\limits^i A, i=0,...,k.$$
In particular, there is a surjection
$$\bigwedge\limits^{c+1}B\to A\otimes\det(C)\to 0,$$
so that, e.g., $A\otimes\det(C)$
is globally generated if
$B$ is.\par
Let $M^d_{\P^n}$ be the vector bundle on $\P^n$ defined
by the natural exact sequence
$$0\to M^d_{\P^n}\to
\L_d\otimes\O_{\P^n}\to\O_{\P^n}(d)\to 0.$$
Then it is well known that $M^1_{\P^n}=\Omega_{\P^n}(1)$
and that $M^1_{\P^n}(1)$ is globally generated.
This can be proved by taking 2nd exterior
power of the Euler sequence. Similarly,
taking $(i+1)$st exterior powers in the above sequence,
we see that
$$(\bigwedge\limits^iM^d_{\P^n})(d)\ \ {\text\it{is\
 globally\ generated}}
\ \ \forall d,i>0.\tag 2.1$$
By surjectivity of the multiplication map
$\L_{d-1}\otimes M^1_{\P^n}\to M^d_{\P^n},$
it follows that
$(\bigwedge\limits^iM^d_{\P^n})(i)$ is globally generated,
hence
$$(\bigwedge\limits^iM^d_{\P^n})(\min(d,i))\ \ {\text\it{is\
 globally\ generated}}
\ \ \forall d,i>0.\tag 2.2$$
More generally, for any Grassmannian
$$G=\Bbb G(r,\P^n),$$
we define a bundle $M^d_G$ by the natural
exact sequence
$$0\to M^d_G\to\L_d\otimes\O_G\to\Sym^d(Q)\to 0,$$
where $Q$ is the tautological
quotient bundle on $G$. Thus $M^1_G$
is just the tautological
subbundle, which fits in an exact sequence
$$0\to M^1_G\to (n+1)\O_G\to Q\to 0.$$
Taking $(n-r+1)$st exterior power of this sequence
and using the fact that $\bigwedge\limits^{n-r}Q=
\det Q=\O_G(1),$ we conclude that
$M^1_G(1)$ is globally
generated, hence it follows as above that
$$M^d_G(1)\ \ {\text\it{is\
 globally\ generated,}} \ \ \forall d\geq 1.
\tag 2.3 $$
\heading 3.\bf{ Local linear systems}\endheading
Fix $x\in\P^n$ and $d$ and set
$$M=M^d_x,$$
i.e. the vector space of homogeneous
polynomials of degree $d$ zero
at $x$. We will sometimes denote
the same vector space by $\L_d(-x)$,
and we similarly have $\L_d(-Z)$ for any subset
$Z\subset \P^n$. A {\it{local linear system at}} $x$ is simply
a subspace $T\subset M$, whose codimension
we denote by $h$.
We define the integer $s=s(T)$ as the smallest
with the property
that, for generic choice of $P_1,...,P_s\in\L_{d-1},$
the natural map
$$\mu_s=\mu_{s,T}:\sum\limits_{i=1}^s P_iM^1_x\to
 M/T\tag 3.1$$
is surjective.
More generally, suppose
$Z$ is an irreducible variety,
$$g:Z\to\P^n$$
a morphism and $T\subset g^*(M^d_{\P^n})$ a subhseaf.
Then we define $s=s(T)$ to be smallest so that the map
$$\mu_s:\sum\limits_{i=1}^s P_ig^*(M^1_{\P^n})
\to g^*(M^d_{\P^n})/T\tag 3.2$$
defined as above is generically surjective
for general choice of $P_i'$s.\par
The interest in studying this integer
comes from the following elementary observation.
\proclaim{Lemma 3.1} In the above situation, we have that
$$H^0(\det (g^*(M^d_{\P^n})/T)\otimes g^*(\O(s)))\neq 0.$$
\endproclaim
\demo{proof} The above sheaf is generically a quotient of
a sum of
tensor products of
sheaves of the form $g^*((\bigwedge\limits^iM^1_{\P^n})(1)).$
\qed\enddemo
Returning to the above situation, set
$$\gamma(i)=\rk(\mu_i)$$
and note that this is a strictly increasing
function for $0\leq i\leq s,$ which is 'concave'
in the sense that $\gamma(i+1)-\gamma(i)$ is non-increasing.
Let $s'=s'_T$ be the smallest $i$ such that
$$\gamma(i+1)-\gamma(i)\leq 1.$$
Note that  $$s'=s\Rightarrow s\leq h/2\tag 3.3$$
while if $s'<s$ we have for generic $P\in\L_{d-1}$
$$\dim(PM^1_x\cap\im(\mu_{s'}))=n-1.\tag 3.4$$
Of course by definition, we have
$$\rk(\mu_{s'})+s-s'=h.\tag 3.5$$
Some interesting things happen when $s-s'>1$:
\proclaim{Lemma 3.2} Suppose that $s-s'>1.$ Then
\item{(i)} There is a line $\ell=\ell_T$ through $x$
such that
$$\L_d(-\ell)=M^d_\ell\subseteq T+\im(\mu_{s'}).$$
\item{(ii)} $\ell$ is uniquely determined independent
of $P_1,...,P_{s'}.$
\item{(iii)} We have
$$\dim(T/(T\cap\L_d(-\ell)))=d-s,\tag 3.6$$
$$\dim(\L_d(-\ell)/(T\cap\L_d(-\ell)))=
\rk(\mu_{s'})-s'.\tag 3.7 $$
\endproclaim
\demo{proof} (i) Let $\ell$ be the unique line through $x$
so that $M^1_\ell\subset M^1_x$ is the codimension-1 subspace
in (3.4).\par
(ii) By assumption, we have for $P,P_1,...,P_{s'+1}$
generic
$$\dim(PM^1_x\cap(T+\sum\limits_1^{s'+1}P_iM^1_x))=
\dim(PM^1_x\cap(T+\sum\limits_1^{s'}P_iM^1_x))=n-1,$$
hence the two intersections are equal and hence
they are independent of $P_{s'+1}$; by symmetry, they are
also independent of any $P_i$, hence so is $\ell.$\par
(iii) Let
$$r=\dim(\frac{M^d_x}{T+M^d_\ell}).$$
Note that the image of $PM^1_x$ in this vector space is
1-dimensional, so if $r\leq s'$ we have that $\mu_{s'}$
is surjective, while if $s'\leq r<s$ we have that
$\mu_r$ is surjective, both of which are impossible.
Since in any case $r\leq s$ we have $r=s$, which
easily implies our first assertion. The second assertion
follows easily by dimension counting.
\enddemo
\heading 4.\bf{ Reduction}\endheading
We are now ready to apply the results of the last section
to our versal family. So let
$$f:   \Y/S'\to\X'=\X\times_SS'$$
be as in \S 1. Then we have the normal sheaf $N_f$ which fits
in an exact sequence
$$0\to T_\Y\to f^*T_{\X'}\to N_f\to 0.$$
By GL$(n+1)$-equivariance, clearly $\Y$ and $\X'$ are
smooth over $\P^n$, so we have another exact
$$0\to T_{\Y/\P^n}\to f^*T_{\X'/\P^n}\to N_f\to 0.$$
Now we have
$$f^*T_{\X'/\P^n}\simeq f^*M^d_{\P^n},$$
therefore we may apply Lemma 3.1
with $T=T_{\Y/\P^n}$
to conclude that
$$H^0(\det(N_f)\otimes f^*\O(s))\neq 0.$$
Now use the generic smoothness of $\Y$ and $\X'$ over $S'$
to conclude that for generic $F\in S'$ we have
$$\det(N_f)|_{Y_F}=\omega_{Y_F}\otimes f^*\omega_{X_F}\inv
=\omega_{Y_F}\otimes f^*\O(-d+n+1).$$
Therefore by Lemma 3.1 we have
$$H^0(\omega_{Y_F}(-d+n+1+s)):
=H^0(\omega_{Y_F}\otimes f^*\O(-d+n+1+s))\neq 0.\tag 4.1$$
The case where $s-s'>1$ will be analyzed at length below;
for now we just note that, in that case, the uniquely
determined line $\l$ going through $f(y)$ for general
$y\in \Y$ gives rise to a rational mapping
$$g:Y\to G=G(1,\P^n).$$
Putting things together we get a lifting of
$f$ to a map
$$\Y\to \Delta:=\{(\l,x,F):x\in\l\cap X_F\}\subset
G\times\P^n\times S'.$$
Viewing $T_{\Y/\P^n}$ as a subsheaf of $\L_d\otimes\O_\Y,$
we have another subsheaf of $\L_d\otimes\O_\Y,$ which we
denote by
$$\L_d(-\l),$$ whose fibre at $y$ is the $d$-forms
vanishing on $\l(y)$;
we set
$$T'=T_{\Y/\P^n}\cap\L_d(-\l).$$
By Lemma 3.1, we have
$$\rk(T_{\Y/\P^n}/T')=d-s\tag 4.2$$
$$\rk(\L_d(-\l)/T')=\rk(\mu_{s'})-s'.\tag 4.3$$
Now suppose that $$s-s'\leq 1.$$
Recall that
$$n-k-1=\rk(\mu_{s'})+s-s'\geq s+s'.$$
Going through the various possibilities we see
easily that in all cases where $s-s'\leq 1$ we have
$$s\leq [\frac{n-k}{2}].\tag 4.4$$
Therefore we conclude, using (4.1),
\proclaim {Alternative 4.1}
In the situation of the main theorem,
we have either
$$H^0(\omega_{Y_F}(-d+n+1+[\frac{n-k}{2}]))\neq 0\tag 4.5$$
or
$$s-s'>1.\tag 4.6$$\endproclaim
Note that when $n-1-k=2,$ the only possibilities are $s'=0,s=2$
and $s=s'=1$. Let us now analyze generally the case $s=1$.
Let $Y$ be a general fibre $Y=Y_F$ and $N$ the normal sheaf
for the map $Y\to\P^n,$ so we have an exact sequence
$$0\to N_f\to N\to \O_Y(d)\to 0.$$
Then choosing a general polynomial $P\in\L_{d-1},$ we get
a diagram with exact rows
$$\matrix 0&\to& f^*M^1_{\P^n}&\to
&\L_1\otimes\O_Y&\to&\O_Y(1)&\to&0\\
           &&\downarrow&&P\downarrow&&\downarrow&&\\
           0&\to&f^*M^d_{\P^n}&\to&\L^d\otimes\O_Y
           &\to&\O_Y(d)&\to&0\\
           &&\downarrow&&\downarrow&&\downarrow&&\\
           0&\to& N_f&\to& N&\to& \O_Y(d)&\to& 0
           \endmatrix$$
Now the assumption that
$s=1$ means that the left column composite
is generically surjective, hence the same is true for the
middle column. Thus we get a generically surjective
map
$$\psi:\L_1\otimes\O_Y\to N$$
which drops rank on Zeroes$(P)$ and also
clearly factors through $H^0(\O_Y(1))\otimes\O_Y$
whose rank we write as $p+1$. Set
$$N'=\im(\psi)/\text{(torsion)}.$$
From the bottom row in the above display we conclude
$$c_1(N')\leq K_Y+(n+2-d)H$$
where $A\leq B$ means $B-A$ is effective and $H$ is
a hyperplane. Now the torsion-free quotient $N'$ of
the trivial bundle corresponds to a rational map to
 a Grassmannian, which
by blowing up we may assume is a morphism
$$\gamma:Y\to G:= G(k+1,p+1),$$
such that
$$\bigcup\limits_{y\in Y}\gamma (y)\text{ spans}\ \C^{p+1}$$
(because $\gamma(y)\ni f(y)$)
and that
$$c_1(N')=\gamma^*(\O_G(1)).$$
 Thus clearly $c_1(N')$ is
effective and if $p>k$ (i.e. $Y$
itself is not a $\P^k$),
then $\gamma$ is nonconstant, so
$$h^0(c_1(N'))\geq 2.$$
It is elementary and well known that
any linear $\P^1$ in $G$
consists of the pencil subspaces contained
in a fixed $k+2$-dimensional subspace
$A\subseteq \C ^{p+1}$ and containing
a fixed $k-$dimensional one $B\subset A$.
Now if $k<p-1$ then $A\subsetneq \C^{p+1}$
and since  $\bigcup\limits_{y\in Y}\gamma (y)$ spans $\C^{p+1}$
it follows that $\gamma (Y)$
cannot be contained in such a pencil, so
$$h^0(c_1(N'))\geq 3.$$
Thus we conclude
\proclaim{Lemma 4.2}
Assuming $s=1$ and that the general $f(Y)$ spans
a $\P^{p}$, we have
$$h^0(\omega_Y(n+2-d))\geq \min(p-1,3).\tag 4.7$$
\endproclaim

Now suppose moreover that $n=4, d=6, k=1$
(still assuming  $s=1$ ). We claim
that if $Y$ is planar, i.e. $p=2$,
then $g(Y)\geq 2$ (in fact, we will show $g(Y)\geq 4$). Let
$$I\subset\L_6\times G(2,\P^4)$$ denote
the (open) set of pairs $(X,B)$ where $X$ is a smooth sextic
hypersurface in $\P^4$
and $B$ is a plane in $\P^4,$ and let $J$ denote the set
of pairs $(D,B)$ where $B$ is a plane in $\P^4$ and $D$
is a sextic curve in $B$, not necessarily reduced or irreducible.
Since a smooth sextic cannot contain a plane, there is a natural
morphism
$$\pi:I\to J,$$
$$(X,B)\mapsto (X\cap B,B).$$
Clearly $\pi$ is a fibre bundle. Now a fundamental fact
of plane geometry  [AC] is that the family of reduced irreducible
plane curves of geometric genus $g$ and degree $d$
is of dimension
$3d+g-1.$ It follows easily from this that the locus
$J_0\subset J$
consisting of pairs $(D,B)$ such that $D$ is the target of
a nonconstant
map from a curve of genus $\leq 3$ is of codimension
$>6$, hence $I_0:=\pi^{-1}(J_0)\subset I$
is also of codimension $>6$.
Since $G(2,\P^4)$ is 6-dimensional,
$I_0$ cannot dominate $\L_6$. Thus we conclude
\proclaim{Lemma 4.3} If $s=1, n=4, d=6, k=1,$ then $Y$
is a curve of genus
at least 3,
except possibly if $Y$ is  genus-2 curve spanning
a hyperplane.\endproclaim


\heading\bf{ 5. Line osculation hierarchy}\endheading
Before continuing with the proof of the main theorem,
we digress to discuss canonical bundles of 'line osculation'
varieties, i.e. the varieties
$$\Delta_r:=\{ (\l,x,F): \l.X_F\geq r.x\}\subset
 G\times\P^n\times S$$
where $G=G(1,\P^n).$ We begin with a formal definition. Let
$Q$ be the tautological quotient bundle on $G$ and
$$\O_S(-1)\subset\L_d\otimes\O_S$$
the tautological subbundle on $S$ (whose fibre at $F$ is $\C F$);
of course, this bundle admits a
canonical nowhere vanishing section, so it is
canonically trivial.
Then $\Delta=\Delta_1$ is the common zero-locus on
$G\times\P^n\times S$ of the  two natural maps
$$\O_{\P^n}(-1)\to Q\tag 5.1$$
$$\O_S(-1)\to \O_{\P^n}(d)\tag 5.2$$
where to save notation we have suppressed the various pullbacks.
From this it is easy to see by the
adjunction formula that $\Delta$,
which is obviously smooth, has canonical bundle
$$\omega_\Delta=\omega_{G\times\P^n\times S}
\otimes \det(Q)\otimes \O_{\P^n}(n-1+d)
=\O_G(-n)\otimes \O_{\P^n}(d-n+1).\tag *$$
Let $\SS$ denote the tautological subbundle on $G$.
Then because (5.1) vanishes on $\Delta$, we have on $\Delta$ a
natural injection
$$\O_{\P^n}(-1)\to \SS$$
whose (rank-1)
quotient we denote by $R^{\text{v}}$. Thus we have an exact
sequence  on $\Delta$
$$0\to R\to\SD\to \O_{\P^n}(1)\to 0\tag 5.3$$
which yields
$$R=\det(\SD)\otimes\O_{\P^n}(-1)
=\O_G(1)\otimes\O_{\P^n}(-1).
\tag 5.4$$
The sequence (5.3) induces a descending filtration $F^.$ on
$\Sym^d(\SD)$ with quotients
$$F^i/F^{i+1}=\O_{\P^n}(d-i)\otimes R^i
=\O_G(i)\otimes\O_{\P^n}(d-2i), 0\leq i\leq d,$$
$$F^{d+1}=0.$$
Note that the vanishing of 5.2 implies that the
natural map
$$\L_d\to\Sym^d(\SS^{\text{v}})$$
induces a map
$$\O_S(-1)\to F^1.$$
We define $\Delta_r$ as the zero scheme of the induced
map
$$\O_S(-1)\to F^1/F^r.$$
In particular, $\Delta_{d+1}$, the zero-scheme of
the natural map
$$\O_S(-1)\to F^1,$$ parametrizes the triples
$$(\l,x,F)$$ such that
$$x\in\l\subset X_F.$$
By considering the projection to (the incidence subvariety
of) $G\times\P^n$, it is easy to see that
$\Delta_r$ is smooth
for $r\leq d+1$. Using the adjunction formula again, we
see that
$$\omega_{\Delta_r}=\omega_\Delta\otimes\det (F^1/F^r)=
\omega_\Delta\otimes\bigotimes\limits_{i=1}^{r-1}
(\O_{\P^n}(d-2i)\otimes \O_G(i)).$$
By (*) we get the formula
$$\omega_{\Delta_r}=\O_G(\frac{r(r-1)}{2}-n)\otimes
\O_{\P^n}(r(d-r+1)-n+1)\tag 5.5$$\par

\heading \bf{6. The case $s-s'>1$}\endheading
We will assume from now on, for the remainder of this paper,
that
$$s'+1<s$$
and that $$H^0(\omega_{Y_F})=0\tag 6.1$$
for generic $F\in S$.
With no loss of generality we may also assume that
$$s>[\frac{n-k}{2}],\tag 6.2$$
and
$$-d+n+1+s>0,\tag 6.3$$
else the Main Theorem's assertion follows from (4.1). This
means
$$2s>n-k, d-n-s-1<0.\tag 6.4$$
Our aim is to show in this case that
$$Y=Y_F$$
is contained in the locus ruled by lines in $X_F$.
We fix a point $x_0\in\P^n$ which
we assume has
coordinates $[1,0,...,0]$, and
let $$G_1 < GL(n+1)$$ be
the stabilizer  of $x_0$, that is,
the parabolic subgroup
$$\left\{ \left[\matrix *&*&\ldots&*\\
0&*&\ldots&*\\
\vdots&\vdots&\vdots&\vdots\\
0&*&\ldots&*\endmatrix\right ]\right\}$$
and let
$$G_0\simeq GL(n)<G_1$$
be the subgroup fixing $X_0$, i.e.
$$G_0=\left\{ \left[\matrix 1&0&\ldots&0\\
0&*&\ldots&*\\
\vdots&\vdots&\vdots&\vdots\\
0&*&\ldots&*\endmatrix\right ]\right\}.$$

Note that the fibre
$Y_0$ of $\Y$ over $x_0$ is invariant under the natural action
of $G_1$  on the
space of polynomials.

Now a homogeneous polynomial of degree $d$
in $X_0,...,X_n$ may be written in the form
$$F=\sum\limits_{|A|\leq d}a_AX_0^{d-|A|}X^A\tag 6.5$$
where
$$A=(A(1),...,A(n))$$
is a multi-index referring
to the variables $X_1,...,X_n.$ Thus the
$X_0^{d-|A|}X^A$ form a basis of the space $\L_d$
of such polynomials and the $a_A$ is the associated
coordinate system or dual basis and give rise to
differentiation operators
$$\p_A=\p/\p a_A.\tag 6.6$$

Note that the action of GL$(n+1)$ on $\Y$ induces an
action of $G_1$ on the fibre $Y_0$ of $\Y$ over
$$x_0=[1,0,...,0]\in\P^n.$$
Working in an analytic neighborhood $Y_{00}$
of a general
point
$$y_0\in Y_0$$
we may identify it via the map $f$ with
an open subset of $S$.
Of course, $S$ being an open subset of $\L_d,$
we have trivializations
$$T_S\simeq\L_d\otimes\O_S, T_S\otimes\O_{Y_0}
\simeq \L_d\otimes \O_{Y_0},\tag 6.7$$
with a frame being provided by the commuting
vector fields $\p_A$ for $|A|\leq d,$
which correspond  to the respective
monomials $X_0^{d-|A|}X^A.$ Then we get an
embedding
$$T_{Y_0}\subset \L_d\otimes\O_{Y_0}.$$
Recall the subbundle
$$T'=T_{Y_0}\cap \L_d(-\l).$$
  Let
$$\L_{d, m}(\-\l)\subset\L_d(\-l)$$
denote the subsheaf of elements of degree
$m$ in $X_1,...,X_n$.
We assume that for $y\in Y_{00}$,
the line $\ell(y)$ is given by
$$X_r=b_r(y)X_1, r=2,...,n,$$
and that at our 'initial' (general) point $y_0\in Y_0$
we have $$b_2(y_0)=...=b_n(y_0)=0.$$

Note that generators for $\L_d(-\ell)$
are given under the identification (6.7)
by the vector fields
$$\p_{Ar}-b_r\p_{A1}\in T_S\otimes\O_{Y_{00}},
|A|\leq d-1, r=2,...n,\tag 6.8$$
corresponding to the polynomials
$$X_0^{d-1-|A|}X^A(X_r-b_rX_1).$$
\proclaim{Lemma 6.1} $T'$ annihilates $b_r$ for all
$r\geq 2.$\endproclaim
\demo{proof}
It suffices to show this at the point $y_0.$
Set
$$T''= T'\cap\L_d(-2\l).$$
We begin by showing that $T''$ annihilates the $b_r.$
Note that sections of $T''$ can be written as
linear combinations of the vector fields
$$\p_{Ars}-b_r\p_{A1s}-b_s\p_{A1r}+b_rb_s\p_{A11},
|A|\leq d-2, r,s\geq 2,$$
 from which it follows
that
$$[T'',T'']\subseteq T'.$$
We thus have a well-defined bracket pairing
$$T''\times(T'/T'')\to T/T'\tag 6.9$$
and we are claiming this vanishes.
$T'/T''$ is generated by sections the form
$$\tau'=\sum\limits_{r=2}^n\sum\limits
_{\alpha=0}^{d-1}c_{r\alpha}(\tau')(\p_{r1^\alpha}-
b_r\p_{1^{\alpha+1}})$$
and the assignment
$$\tau'\mapsto (c_{r\alpha}(\tau'))$$
identifies $T'/T''$ with a subsheaf, which we may assume
is locally free,
$$T'/T''\subseteq M_{(n-1)\times d}(\O)$$
of $\O_{Y_{00}}-$ valued  $(n-1)\times d$ matrices, and the
corank of this subsheaf is at most the
corank of $T'$ in $\L_d(-\l)$.
Hence by (4.3) this corank is
at most $n-k-1-s$. On the other hand a section
$$\tau\in T/T'$$ can be written in the form
$$\tau=\sum\limits_{\beta=1}^d c_\beta(\tau)
\p_{1^\beta}$$
and this identifies $T/T'$ with a subsheaf
$$T/T'\subseteq d\O_{Y_{00}}$$
whose corank by (4.2) is exactly $s$. Now for $\tau'\in T'/T'',
\tau''\in T'',$ the pairing (6.9)
 is given by
$$[\tau',\tau'']=
\sum c_{r\alpha}\tau''(b_r)\p_{1^{\alpha+1}}.$$
Now suppose $\tau''(b_r)\neq 0$ for some $r$.
Then the map
$$M_{(n-1)\times d}(\O)\to d\O$$
given by multiplication by the vector
$$\beta(\tau''):=(\tau''(b_2),...,\tau''(b_r))$$
is surjective, hence its restriction
$$T'/T''\to T/T'$$
can only lower corank, i.e. has image of
corank at most $$n-k-1-s$$ in $d\O$.
But this contradicts the
fact that $T/T'$ has corank $s$ in $d\O$,
while by (6.2) we have $s>n-k-1-s$.
\par
This proves that
$T''$ annihilates all the $b_r,$ which in turn implies that
we have a well-defined bracket pairing
$$T'/T''\times T'/T''\to T/T',$$
which has the form
$$[\tau_1,\tau_2]=\sum\limits_{s\beta}\tau_1(b_s)
c_{s\beta}(\tau_2)\p_{1^{\beta+1}}-
\sum\limits_{r\alpha}\tau_2(b_s)c_{r\alpha}(\tau_1)
\p_{1^{\alpha+1}}.\tag 6.10$$
Now suppose that $$\tau_1(b_s)\neq 0$$
for some
$$\tau_1\in T'/T'', s\geq 2.$$
Then the following linear algebra
observation (copmare with [P]) shows that the set of
elements $$[\tau_1,\tau_2]$$ ranges over a subspace
of $T/T'$ of codimension at  most
$$1+n-k-1-s,$$
which by (6.2) is $>s$, so we have a contradiction as above.
\qed

\enddemo
\proclaim{Sublemma 6.1.1} Let
$$U,V$$ be finite-dimensional
vector spaces, $$H<{\text{\rm{Hom}}}(U,V)$$
a subspace of codimension
$c$, and
$$\beta:H\to U$$ a nonzero homomorphism. Define
a pairing
$$I(\beta):\bigwedge\limits^2H\to V$$
by
$$A\wedge B\mapsto A\phi(B)-B\phi(A).\tag 6.11$$
Then the image of $I(\beta)$ is of codimension at most
$c+1$ in $V$.\endproclaim
\demo{proof} Pick $c$ rank-1 elements in Hom$(U,V)$
that are linearly independent modulo $H$,
and extend $\beta$ linearly to be zero on these.
It is clear that this extension increases the rank of
$I(\beta)$ by at most $c$. Therefore it suffices to
prove the Sublemma in case $c=0$. In this case,
pick any $A$ of rank 1 with $\beta(A)\neq 0$. Then
as $B$ ranges over $H={\text{Hom(U,V)}},$ the first term in
(6.11) ranges over a 1-dimensional subspace while
the second ranges over all of $V$, hence the difference
ranges at least over a codimension-1 subspace.
\qed\enddemo
As an immediate consequence of Lemma 6.1, we conclude
\proclaim{Lemma 6.2} The distribution $T'$ is integrable
\endproclaim\demo{proof}
For $\tau_1, \tau_2\in T'$, write
$$\tau_i=\sum c_{Ar}(\tau_i)(\p_{Ar}-b_r\p_{A1}),
i=1,2.$$
Then we calculate
$$[\tau_1,\tau_2]\equiv\sum (c_{r1^\alpha}(\tau_1)
\tau_2(b_r)-c_{r1^\alpha}(\tau_2)\tau_1(b_r))
\p_{1^{\alpha+1}}\mod T'$$
and by Lemma 6.1 the latter vanishes.\qed\enddemo

\proclaim{Lemma 6.3}
 For general $y\in Y,$ the line $\ell(y)$
is either contained in
the hypersurface $X_{F(y)}$ or meets it
set-theoretically in
at most 2 points.

\endproclaim
\demo{proof}
emo{proof}
The natural map
$$\l:Y_0\to G$$
is clearly $G_1-$equivariant, hence
generically onto the locus of lines through
$x_0$. Let $Z$ be its fibre over
 $$\ell_0=\{X_2=...=X_n=0\},$$
 which we may assume is a general fibre, hence of
 codimension $n-1$ in $Y_0$. There is a natural
 map
 $$R: Z\to H^0(\O_{\l_0}(d)(-x_0)),$$
 $$R(F)=F|_{\l_0}.$$
 Clearly $R$ is equivariant with respect to the stabilizer
 $G_{11}$ of $(x_0,\l_0)$ in GL$(n)$, so its image is
 a cone invariant under the stabilizer $H$ of $x_0$
 in GL$(\l_0)$. Note that our distribution $T'$
 is none other than the vertical tangent space for
 the map $R$. From (4.2) it follows that
 $$\dim R(Z) =d-s-n+1$$
 and by our assumptions (6.3) this is $\leq 1$.
 As $H$ acts doubly transitively on $\l -x_0$, it follows
 easily that $R(Z)$ is contained in the set of
 multiples of $X_1^d$.
 Therefore the lemma holds.\qed
 \remark{Remark 6.3.1} Note that if we assume instead of (6.3)
 the stronger condition that
 $$-d+n+1+s>1,$$
 then we conclude in the same way that
$$\l(y)\subset X_{F(y)}. $$
In the next Lemma we will see that this holds anyway,
without the extra hypothesis.\par
On the other hand if we assume the weaker condition that
$$-d+n+1+s>-1,\tag *$$
then the same argument, using double transitivity of $H$
shows that $\l(y)$ meets $X_{F(y)}$ in at most 2 points
set-theoretically, hence is either contained in $X_{F(y)}$,
or is $d-$fold tangent to it at $x(y)$, or has contact of order
$r$ with it at $x(y)$ and of order $d-r$ at another point
$x'$, for some $1\leq r\leq d-1.$
Note that, as above, the inequality (*) may
be assumed to hold if $$h^0(Y_F(-1))=0.$$
\endremark
\enddemo
\
\proclaim{Lemma 6.4 }
 For general $y\in Y,$ the line $\ell(y)$
is contained in the hypersurface $X_{F(y)}.$\endproclaim
\demo{proof} Denote by $$\pi:Y\to I$$ the natural map of $Y$
to the incidence variety of pairs (point on line), which is
clearly generically submersive (e.g. by $GL(n+1)-$
equivariance) and let
$T_\pi$ denote the vertical tangent space for this map.
Note that $$F\in T_\pi$$ (at $(x,l,F)$): indeed this is immediate
from the fact that our universal hypersurface $X/S$ together,
we may assume, with its subvariety $Y/S$ come from analogous
families defined over the projectivization $\P(S)$ so we
may assume the map $\pi$ descends accordingly; on the other hand
$F$ as tangent vector dies in the map $S\to \P(S).$\par
Now let $$M$$ be the vector bundle on $Y$ whose fibre at
$$y\mapsto (x,\ell,F)$$ is
$$\L(-\ell),$$ the polynomials of degree $d$
(on $\P^n$) vanishing on $\ell$, considered as a subbundle
of the tangent bundle to $S$ (pulled back to $Y$). Suppose
first that
$T_\pi\cap M$ has generic corank $\leq c-1$ in $M$,
where $$c:=2n-1-d-k.$$ Then the rank of $T_\pi\cap M$
is at least
$$N-(d+c)=N+k-(2n-1).$$
But owing to the generic submersivity of $\pi$ the latter is
precisely the rank of $T_\pi$. Hence $$T_\pi\subseteq M,$$
and therefore $F\in M,$ i.e. $F$ vanishes on $\ell$ as
claimed.\par
Now suppose that $T_\pi\cap M$ has (generic) corank $c$
in $M$, which clearly implies that the natural map
$$M\to N$$
of $M$ to the normal bundle of $Y$ in $\Delta_d$ is
generically surjective, hence so is $$\bigwedge\limits^c M
\to\bigwedge\limits^cN=\det N.$$
Now assume that the line $\l$ has $d-$fold contact if $X_F$
at $x_0$.
Since
we have $d>n, d(d-1)/2-n\geq c$,
the results of \S 2 and \S 5,(5.5) imply that
$\omega_{\Delta_d}\otimes
\bigwedge\limits^cM$ is globally generated and
it follows that
 $$\omega_{\Delta_d}\otimes\O_{\P^n}(-1)
 \otimes\det N\subseteq \omega_Y(2-d)$$
has a nonzero section , hence so does
$\omega_{Y_F}(-1)$ for generic $F$, which
is a contradiction.\par
\qed\enddemo
\heading 7. The case $s'=0$\endheading
Throughout this section we continue to assume
the hypotheses of the previous section,
and assume additionally that $$s'=0.$$
As has been remarked before, the latter assumption holds
automatically whenever $$s>1, k= n-3.$$
Our aim is to show,
under these assumptions, that a general $Y_F$ is
ruled by lines. This will complete the proof
of the Main Theorem.\par
Note that the condition $s'=0$ is equivalent to
$$T'=\L_d(-\l),$$
i.e. to the assumption that
$$\L_d(-\l)\subseteq T_{\Y/\P^n}.\tag 7.1$$
\proclaim{Lemma 7.1} $\ell(y)$ is tangent to $Y_F$ for general
$y\in Y_F.$ \endproclaim
\demo{proof} Using notations as in the previous section,
let
$$v\in gl(n+1)$$ be
the vector field
$$v=X_0\p/\p X_1,$$
which at $y_0$ is nonzero and points in the direction of
$\ell_0,$
and set $$G=v(F).$$ Considering $G$ as an element of
$T_FS$, note that by $GL(n+1)-$equivariance,
the normal field $$\bar{v}$$ to $Y=Y_F$ corresponding
to $G$ is just the natural image of $v$. On
the other hand since $$G|_\ell=0,$$
our assumption $s'=0,$ in the form (7.1),
yields that $(G,0)$ is tangent to
$\Y$, hence $$\bar{v}(y_0)=0.$$
 Thus $v$ is tangent to $Y_F$
$y_0,$ hence so is $\ell_0$, as claimed.\qed\enddemo
\proclaim{Lemma 7.2} $Y_F$ is ruled by lines.\endproclaim
\demo{proof} In the above notations, note that since the
action of $v$ on the Grassmannian fixes $\ell_0$, it follows
by equivariance that on $\Y,$
$$\ell=\ell(y)$$ is fixed to first order
as $y$
 moves out of $y_0$ in the direction
$$v(y_0)\in T_{y_0}\Y.$$ Viewing the latter
as a subspace of $T_FS\times T_{y_0}Y_F,$ write
$$v(y_0)=(G,v_Y),$$
where $$G=v(F)\in\L_d, v_Y\in T_{y_0}Y_F.$$
As we have seen, $G$ is in fact in
$\L_d(-\l)$ which under our assumption
coincides with $T'(y_0)$ (more precisely,
$(G,0)\in T'(y_0)$). Therefore by Lemma 6.1, $\l$
is constant to 1st order in the direction
$(G,0).$ Therefore it is also constant to 1st
order in the direction $(0,v_Y)$.
Now let
$$u$$ be a vector
field locally on $Y_F$ whose direction
at any $y$ is the same as $\ell(y)$.
Then along any integral
arc $U$ of $u$,
$\ell$ is infinitesimally constant, hence
constant,
hence $U$ itself in on a line. Thus $Y_F$ is ruled by
lines.\qed\enddemo
{\bf {Acknowledgement}}\
The first author thanks the Institute for Advanced Study for its
hospitality during the period during which his
contributions this work were
completed.
The second author's work
on these questions was begun while he was visiting
Taiwan under the auspices of the National Center for Theoretical
Sciences, J. Yu, Director. He would like to thank the NCTS for the
invitation, and especially Prof. L.C.Wang and his colleagues at
National Dong-Hwa University for their outstanding hospitailty.


\enddocument